\documentclass[12pt,twoside]{article}
\usepackage{amsmath,amsfonts,amssymb,theorem}
\textheight215truemm \textwidth 138truemm  \date{}
 \newtheorem{thm}{Theorem}[section]
 \newtheorem{cor}[thm]{Corollary}
 \newtheorem{lem}[thm]{Lemma}
 \newtheorem{prop}[thm]{Proposition}

 \theorembodyfont{\rmfamily}
 \newtheorem{rem}[thm]{Remark}
 \newtheorem{spec}[thm]{Speculation}

 \newenvironment{pf}{\paragraph{Proof}}{\par\medskip}

 \newcommand{\wave}{\widetilde}
 \newcommand{\iso}{\cong}
 \newcommand{\rest}[1]{_{{\textstyle{|}}#1}} 
 \newcommand{\rA}[1]{\mathrm A_{#1}}
 \newcommand{\rcA}[1]{\mathrm{cA}_{#1}}
 \newcommand{\rD}[1]{\mathrm D_{#1}}
 \newcommand{\rcD}[1]{\mathrm{cD}_{#1}}
 \newcommand{\rcE}[1]{\mathrm{cE}_{#1}}

 \newcommand{\C}{\mathbb C}
 \newcommand{\F}{\mathbb F} 
 \newcommand{\PP}{\mathbb P}
 \newcommand{\R}{\mathbb R}
 \newcommand{\Z}{\mathbb Z}

 \newcommand{\ind}{\operatorname{ind}} 
 \newcommand{\Pic}{\operatorname{Pic}}
 \newcommand{\Supp}{\operatorname{Supp}}
 \newcommand{\Tors}{\operatorname{Tors}}
 \newcommand{\coker}{\operatorname{coker}}

 \newcommand{\HR}[2]{H^{#1}(#2,\R)}
 \newcommand{\HZ}[2]{H^{#1}(#2,\Z)}

 \newcommand{\Kbar}{\overline{\mathcal K}}

 \newcommand{\Bbar}{\overline B}
 \newcommand{\Ubar}{\overline U}
 \newcommand{\Xbar}{\overline X}
 \newcommand{\delbar}{\overline\partial}
 \newcommand{\fiemapX}{\fie\colon X\to\Xbar}

 \newcommand{\cD}{\mathcal D}
 \newcommand{\cE}{\mathcal E}
 \newcommand{\cF}{\mathcal F}
 \newcommand{\cH}{\mathcal H}
 \newcommand{\cM}{\mathcal M}
 \newcommand{\cMbar}{\overline{\mathcal M}}
 \newcommand{\Oh}{\mathcal O}
 \newcommand{\cS}{\mathcal S}
 \newcommand{\cU}{\mathcal U}
 \newcommand{\cW}{\mathcal W}
 \newcommand{\cX}{\mathcal X}
 \newcommand{\cY}{\mathcal Y}
 \newcommand{\cZ}{\mathcal Z}
 \newcommand{\cK}{\mathcal K}

 \newcommand{\Si}{\Sigma}
 \newcommand{\Om}{\Omega}
 \newcommand{\De}{\Delta}
 \newcommand{\fie}{\varphi}
 \newcommand{\ep}{\varepsilon}

 \author{P.M.H. Wilson}
 \title{Flops, Type~III contractions and Gromov--Witten invariants\\ on
Calabi--Yau threefolds}

\begin{document}
 \maketitle

 \section{Introduction}

In this paper, we investigate Gromov--Witten invariants associated to
exceptional classes for primitive birational contractions on a Calabi--Yau
threefold $X$. As already remarked in \cite{18}, these invariants are
locally defined, in that they can be calculated from knowledge of an open
neighbourhood of the exceptional locus of the contraction; intuitively, they are 
the numbers of rational curves in such a neighbourhood. In \S\ref{sec1},
we make this explicit in the case of Type~I contractions, where the
exceptional locus is by definition a finite set of rational curves.
Associated to the contraction, we have a flop; we deduce furthermore in
Proposition~\ref{prop_1.1} that the changes to the basic invariants (the
cubic form on $H^2(X,\Z)$ given by cup product, and the linear form given by
cup product with the second Chern class $c_2$) under the flop are explicitly
determined by the Gromov--Witten invariants associated to the exceptional
classes.

The main results of this paper concern the Gromov--Witten invariants
associated to classes of curves contracted under a Type~III primitive
contraction. Recall \cite{17} that a primitive contraction $\fiemapX$ is of
Type~III if it contracts down an irreducible divisor $E$ to a curve of
singularities $C$. For $X$ a smooth Calabi--Yau threefold, such contractions
were studied in \cite{18}; in particular, it was shown there that the curve
$C$ is smooth and that $E$ is a conic bundle over $C$. We denote by
$2\eta\in H_2(X,\Z)/\Tors$ the numerical class of a fibre of $E$ over $C$. In
the case when $E$ is a $\PP ^1$-bundle over $C$, this may in fact be a
primitive class, and so the notation is at slight variance with that adopted
in \S\ref{sec1}, where $\eta$ is assumed to be the primitive class. In the
case when the class of a fibre is not primitive (for instance, when $E$ is not
a $\PP ^1$-bundle over $C$), the primitive class contracted by $\fie$ will be
$\eta$. We denote the Gromov--Witten numbers associated to $\eta$ and
$2\eta$ by $n_1$ and $n_2$, with the convention that $n_1=0$ if $2\eta$ is the
primitive class. The above conventions have been adopted so as to achieve
consistency of notation for all Type~III contractions.

If the genus $g$ of the curve $C$ is strictly positive, under a
general holomorphic deformation of the complex structure on $X$, the
divisor $E$ disappears leaving only finitely many of its fibres, and (except
in the case of elliptic quasiruled surfaces, where all the Gromov--Witten
invariants vanish) we have a Type~I contraction. The results of \S\ref{sec1}
may then be applied to deduce the Gromov--Witten invariants associated to
the classes $m\eta$ for $m>0$. These are all determined by the Gromov--Witten
numbers $n_1$ and $n_2$, and explicit
formulas for $n_1$ and $n_2$ are given in Proposition~\ref{prop_2.3}; in
particular $n_2=2g-2$.

The formulas for $n_1$ and $n_2$ remain valid also for $g=0$, although the
slick proof given in Proposition~\ref{prop_2.3} for the case $g>0$ no longer
works. The formula for $n_1$ is proved for all values of $g(C)$ by local
deformation arguments in Theorem~\ref{thm_2.5}. Verifying that $n_2=-2$ in
the case when $g(C)=0$ is rather more difficult, and involves the technical
machinery of moduli spaces of stable pseudo\-holomorphic maps and the virtual
neighbourhood method, as used in \cite{2,9} in order to construct
Gromov--Witten invariants for general symplectic manifolds. In particular,
we shall need a cobordism result from \cite{13}, which we show in
Theorem~\ref{thm_3.1} applies directly in the case where no singular fibre
of $E$ is a double line. The general case may be reduced to this one by
making a suitable almost complex small deformation of complex structure. In
\S\ref{sec4}, we give an application of our calculations. In \cite{18}, it
was shown that if $X_1$, $X_2$ are Calabi--Yau threefolds which are
symplectic deformations of each other (and general in their complex moduli),
then their K\"ahler cones are the same. Now we can deduce
(Corollary~\ref{cor4.1}) that corresponding codimension one faces of these
cones have the same contraction type.

The author thanks Yongbin Ruan for the benefit of conversations concerning
material in \S\ref{sec3} and his preprint \cite{13}. 

 \section{Flops and Gromov--Witten invariants}\label{sec1}

If $X$ is a smooth Calabi--Yau threefold with K\"ahler cone $\cK$, then the
nef cone $\Kbar$ is locally rational polyhedral away from the cubic cone
 \[
W^*=\bigl\{ D\in \HR 2 X \ ; \ D^3=0\bigr\};
 \]
moreover, the codimension one faces of $\Kbar$ (not contained in $W^*$)
correspond to primitive birational contractions $\fiemapX$ of one of three
different types \cite{17}.

In the numbering of \cite{17}, Type~I contractions are those where only a
finite number of curves (in fact $\PP^1$s) are contracted. The singular
threefold $\Xbar$ then has a finite number of cDV singularities. Whenever one
has such a small contraction on $X$, there is a flop of $X$ to a different
birational model $X'$, also admitting a birational contraction to $\Xbar$;
moreover, identifying $\HR 2 {X'}$ with $\HR 2 {X}$, the nef cone of $X'$
intersects the nef cone of $X$ along the codimension one face which defines
the contraction to $\Xbar$ \cite{6, 7}. It is well known \cite{7} that
$X'$ is smooth, projective and has the same Hodge numbers as $X$, but that
the finer invariants, such as the cubic form on $\HZ 2 X$ given by cup
product, and the linear form on $\HZ 2 X$ given by cup product with
$c_2(X)=p_1(X)$, will in general change. Recall that, when $X$ is simply
connected, these two forms along with $\HZ 3 X$ determine the diffeomorphism
class of $X$ up to finitely many possibilities \cite{14}, and that if
furthermore $H_2(X,\Z)$ is torsion free, this information determines the
diffeomorphism class precisely \cite{16}.

When the contraction $\fiemapX$, corresponding to such a {\em flopping face}
of $\Kbar$, contracts only isolated $\PP^1$s with normal bundle $(-1,-1)$
(that is, $\Xbar$ has only simple nodes as singularities), then it is a
standard calculation to see how the above cubic and linear forms (namely the
cup product $\mu\colon \HZ 2 X\to\Z$, and the form $c_2\colon \HZ 2 X\to\Z$)
change on passing to $X'$ under the flop. Since any flop is an isomorphism
in codimension one, we have natural identifications 
 \[
 \HR 2 {X'}\iso \Pic_{\R}(X')\iso \Pic_{\R}(X)\iso \HR 2 X.
 \]
If
we are in the case where the exceptional curves $C_1,\dots,C_N$ are
isolated $\PP^1$s with normal bundle $(-1,-1)$, and if we denote by $D'$ the
divisor on $X'$ corresponding to $D$ on $X$, then
 \[
 (D')^3=D^3-\sum(D\cdot C_i)^3
\quad \text{and}\quad c_2(X')\cdot D'=c_2(X)\cdot D+2\sum D\cdot C_i \ .
 \]
 This is an easy verification -- see for instance \cite{1}.

 \begin{prop}\label{prop_1.1} Suppose that $X$ is a smooth Calabi--Yau
threefold, and $\fiemapX$ is any Type~I contraction, with $X'$ denoting the
flopped Calabi--Yau threefold. The cubic and linear forms $(D')^3$ and
$D'\cdot c_2(X')$ on $X'$ are then explicitly determined by the cubic and
linear forms $D^3$ and $D\cdot c_2(X)$ on $X$, and the $3$-point
Gromov--Witten invariants $\Phi_A $ on $X$, for $A\in H_2(X,\Z)$ ranging
over classes which vanish on the flopping face.
 \end{prop}

 \begin{rem} This is essentially the statement from physics that the A-model
3-point correlation function on $\cK(X)$ may be analytically continued to
give the A-model 3-point correlation function on $\cK(X')$.
 \end{rem}

 \begin{pf} We use the ideas from \cite{18}; in particular, we know that on
a suitable open neighbourhood of the exceptional locus of $\fie$, there
exists a small holomorphic deformation of the complex structure for which
the exceptional locus splits up into disjoint $(-1,-1)$-curves (\cite{18},
Proposition~1.1).

Let $A\in H_2(X,\Z)$ be a class with $\fie_* A=0$. The argument from
\cite{18}, Section~1 then shows how the Gromov--Witten invariants
$\Phi_A(D,D,D)$ can be calculated from local information. Having fixed a
K\"ahler form $\omega$ on $X$, a small deformation of the holomorphic
structure on a neighbourhood of the exceptional locus may be patched
together in a
$C^{\infty}$ way with the original complex structure to yield an almost
complex structure tamed by $\omega$, and the Gromov--Witten invariants can
then be calculated in this almost complex structure. The Gromov Compactness
Theorem is used in this argument to justify the fact that all of the
pseudo\-holomorphic rational curves representing the class $A$ have images
which are $(-1,-1)$-curves in the deformed local holomorphic structure.

Here we also implicitly use the Aspinwall--Morrison formula for the
contribution to Gromov--Witten invariants from multiple covers of
infinitesimally rigid $\PP^1$s, now proved mathematically by Voisin
\cite{15}. So if $n(B)$ denotes the number of $(-1,-1)$-curves representing
a class given $B$, then
 \[
 \Phi_A(D,D,D)=(D\cdot A)^3 \sum_{kB=A} n(B)/k^3,
 \]
where the sum is taken over all integers $k>0$ and classes $B\in H_2(X,\Z)$
such that $kB=A$. So if $H_2(X,\Z)$ is torsion free and $A$ is the primitive
class vanishing on the flopping face, this says that 
 \[
 \Phi_{mA}(D,D,D)=(D\cdot A)^3 \sum_{d|m} n(dA)d^3.
 \]
Recall that the Gromov--Witten invariants used here are the ones (denoted
$\wave\Phi$ in \cite{12}) which count marked parametrized curves
satisfying a perturbed pseudo\-holomorphicity condition. Knowledge of the
numbers $n(A)$ for the classes A with $\fie_* A=0$ determines the
Gromov--Witten invariants $\Phi_{A}$ for classes A with $\fie_*A=0$, and
vice-versa.

If we can now show that the local contributions to $(D')^3$ and $D'\cdot
c_2(X')$ are well-defined and invariant under the holomorphic deformations
of complex structure we have made locally, then the obvious formulas for
them will hold. Let $\eta\in H_2(X,\Z)/\Tors$ be the primitive class with
$\fie_* \eta=0$ and $n_d$ denote the total number of $(-1,-1)$-curves on the
deformation which have numerical class $d\eta$; the $n_d$ are therefore
nonnegative integers (cf.\ \cite{10}, Remark~7.3.6). Then
 \begin{align}
 (D')^3&=D^3-(D\cdot\eta)^3\sum_{d>0}n_dd^3, \tag{2.1.1}\\
 D'\cdot c_2(X')&=D\cdot c_2(X)+2(D\cdot\eta)\sum_{d>0}n_dd. \tag{2.1.2}
 \end{align}
To justify the premise in the first sentence of the paragraph, the basic result
needed is that of local conservation of number, as stated in \cite{3},
Theorem~10.2.

For calculating the change in $D^3$ for instance, let $X$ now denote the
neighbourhood of the exceptional locus of $\fie$ and $\pi\colon \cX\to B$ the
small deformation under which the exceptional locus splits up into
$(-1,-1)$-curves. So we have a regular embedding (of codimension six)
 \[
 \renewcommand{\arraystretch}{1.3}
 \begin{matrix}
 \cX &\hookrightarrow & \cX \times \cX \times \cX &=\cY \\
 \downarrow && \downarrow \\
 B & \kern1.2em=\kern-1.2em & B 
 \end{matrix}
 \]
In order to calculate the triple products $D_1'\cdot D_2'\cdot D_3'$ from 
$D_1\cdot D_2\cdot D_3$ and the numbers $n_d$, we may assume {\em wlog}
that the $D_i$ are very ample, and so in particular we get effective
divisors $\cD_1$, $\cD_2$ and $\cD_3$ on $\cX /B$. Applying \cite{8},
Theorem~11.10, we can flop in the family $\cX\to B$, hence obtaining a
deformation $\cX'\to B$ of the flopped neighbourhood $X'$. We wish to
calculate the local contribution to $D_1'\cdot D_2'\cdot D_3'$; with the
notation as in \cite{3}, Theorem~10.2, we have a fibre square 
 \[
 \renewcommand{\arraystretch}{1.3}
 \begin{matrix}
 \cW & \longrightarrow & \cD_1' \times \cD_2'\times \cD_3' \\
 \downarrow && \downarrow \\
\cX' &\longrightarrow & \cX' \times \cX'\times \cX' 
 \end{matrix}
 \]
with $\Supp(\cW)=\bigcap \Supp(\cD_i')$. Furthermore, we may assume that the
divisors $\cD_i $ were chosen so that $\cD_1\cap\cD_2\cap\cD_3$ has no
points in $\cX$, and so in particular $\cW$ is proper over $B$. Letting
$D_i'(t)$ denote the restriction of $\cD_i'$ to the fibre $X_t'$, we
therefore have a well-defined local contribution to $D_1'(t)\cdot
D_2'(t)\cdot D_3'(t)$ (concentrated on the flopping locus of $X_t'$), which
is moreover independent of $t\in B$. Thus by making the local calculation as
in (7.4) of \cite{1}, we deduce that
 \[
 D_1\cdot D_2\cdot D_3-D_1'\cdot D_2'\cdot D_3'
 =(D_1\cdot \eta)(D_2\cdot \eta)(D_3\cdot \eta)
\sum_{d>0} n_d d^3 
 \]
 as required.

The proof for $c_2\cdot D$ is similar. Here we consider the graph 
$\wave X \subset X \times X'$ of the flop, with $\pi_1\colon \wave X\to X
$ and $\pi_2\colon \wave X\to X' $ denoting the two projections, and 
$E \subset \wave X$ the exceptional divisor for both $\pi_1$ and $\pi_2$. 
Then $\pi_2^*(T_{X'})\rest{\wave X \setminus E}=\pi_1^*(T_X)\rest{\wave X
\setminus E} $, and so in particular $c_2(\pi_2^* T_{X'})-c_2(\pi_1^*T_X)$
is represented by a 1-cycle $Z$ on $E$. Suppose {\em wlog} that $D$ is very
ample, and that $D'$ denotes the corresponding divisor on $X'$. Set $\pi_1^*
D=\wave D$ and $ \pi_2^* D'=\wave D +F$, with $F$ supported on $E$. Then
$c_2(X')\cdot D'=c_2(\pi_2^* T_{X'})\cdot(\wave D+F)$. Hence
 \[
 c_2(X')\cdot D'-c_2(X)\cdot D=c_2(\pi_2^* T_{X'})\cdot F+Z\cdot
\wave D=c_2((\pi_2^* T_{X'})\rest F)+(Z\cdot \wave D)_E 
 \]
where the right-hand side is purely local. Note the slight abuse of notation
here that $F$ denotes also the fixed {\em scheme} for the linear system $
|\pi_2^* D' |$.

Now taking $X$ to be a local neighbourhood of the flopping locus, and taking 
a small deformation $\cX\to B$ as before, we obtain families $\cX'$,
$\wave \cX $, $\cD$, $\cE$, $\cF$ and $\cZ$ over B (corresponding to
$X'$, $\wave X$, $D$, $E$, $F$ and $Z$). For ease of notation, we shall use
$\pi_1$ and $\pi_2$ also for the morphisms of families $\wave \cX\to
\cX$, respectively $\wave \cX\to \cX'$. Applying \cite{3},
Theorem~10.2 to the family of vector bundles $(\pi_2^* T_{\cX'/B})\rest{\cF}$
on the scheme $\cF$ over $B$ yields that $c_2((\pi_2^*
T_{\cX'/B})\rest{F_t})$ is independent of $t\in B$. Noting that $\wave\cD
\hookrightarrow \wave\cX$ is a regular embedding, we apply the same
theorem to the fibre square
 \[
 \renewcommand{\arraystretch}{1.3}
 \begin{matrix}
 \wave\cD \times_{\wave \cX} \cE & \longrightarrow & \cE \\
\downarrow && \downarrow \\
 \wave\cD &\longrightarrow & \wave\cX
 \end{matrix}
 \]
and the cycle $\cZ$ on $\cE$. This yields that $(Z_t\cdot\wave D_t)_{E_t}$
on $E_t$ is independent of $t\in B$, where by definition 
 \[
 Z_t=c_2(\pi_2^* T_{\cX'/B})\rest{X_t}-c_2(\pi_1^* T_{\cX /B})
\rest{X_t}.
 \]
 Thus the local contribution to $D'(t)\cdot c_2(X_t')$ is
well-defined and independent of $t$, and so we need only make the local
calculation for generic $t$ (where the exceptional locus of the flop consists of
disjoint $(-1,-1)$-curves). This calculation may be found in \cite{1}, (7.4).
 \end{pf}

 \begin{spec} There are reasons for believing that only the numbers $n_1$
and $n_2$ are nonzero, and hence that the Gromov--Witten invariants
associated to classes $m\eta$ for $m>2$ all arise from multiple covers. If
this speculation is true, then the changes under flopping to the cubic form
and the linear form would be determined by these two integers, and
conversely.
 \end{spec}

 \section{Type~III contractions and Gromov--Witten invariants}\label{sec2}

The main results of this paper concern the Gromov--Witten invariants
associated to classes of curves contracted under a Type~III primitive
contraction. Recall \cite{17} that a primitive contraction $\fiemapX$ is of
Type~III if it contracts down an irreducible divisor $E$ to a curve of
singularities $C$. For $X$ a smooth Calabi--Yau threefold, such contractions
were studied in \cite{18}; in particular, it was shown there that the curve
$C$ is smooth and that $E$ is a conic bundle over $C$. We denote by
$2\eta\in H_2(X,\Z)/\Tors$ the numerical class of a fibre of $E$ over $C$. As
explained in the Introduction, we denote by $n_1$ and $n_2$ the Gromov--Witten
numbers associated to the classes $\eta$ and $2\eta$, where 
$n_1=0$ if $E$ is a $\PP^1$-bundle over $C$. If the generic fibre of $E$
over $C$ is reducible (consisting of two lines, each with class $\eta$), then,
except in two cases, it follows from the arguments of \cite{18}, \S4 that, by
making a global holomorphic deformation of the complex structure, we may
reduce down to the case where the generic fibre of $E$ over $C$ is
irreducible. The two exceptional cases are:
 \begin{enumerate}
 \renewcommand{\labelenumi}{(\alph{enumi})}
 \item $g(C)=1$ and $E$ has no double fibres.
 \item $g(C)=0$ and $E$ has two double fibres.
 \end{enumerate}

However, Case~(a) is an {\em elliptic quasi-ruled} surface in the terminology
of \cite{18}, and hence disappears completely under a generic global
holomorphic deformation. In particular, we know that all the Gromov--Witten
invariants $\Phi_A$ are zero, for $A\in H_2(X,\Z)$ having numerical class
$m\eta$ for any $m>0$.

In Case~(b), $E$ is a nonnormal generalized 
del Pezzo surface $\overline\F_{3;2}$ of
degree 7 (see \cite{18}). As argued there however, we may make a holomorphic
deformation in a neighbourhood of $E$ so that $E$ deforms to a {\em smooth}
del Pezzo surface of degree 7, and where the class $\eta$ is then represented
by either of two `lines' on the del Pezzo surface (which are $(-1,-1)$-curves
on the threefold); hence $n_1=2$. In fact, the smooth del Pezzo surface is
fibred over $\PP^1$ with one singular (line pair) fibre. The arguments we
give below may be applied locally (more precisely with the global almost
complex stucture obtained by suitably patching the local small holomorphic
deformation on an open neighbourhood of $E$ with the original complex
structure), and the Gromov--Witten invariants may be calculated as if the
original contraction $\fie$ had contracted such a smooth del Pezzo surface
of degree 7. In particular, $n_1=2$ comes from the two components of the
singular fibre (Theorem~\ref{thm_2.5}), and $n_2=-2$ is proved in
\S\ref{sec3} (see also Remark~\ref{rem_2.4}).

Let us therefore assume that the generic fibre of $E$ over $C$ is
irreducible, and so in particular $E\to C$ is obtained from a $\PP^1$-bundle
over $C$ by means of blowups and blowdowns. Moreover $E$ itself is a conic
bundle over $C$, and so its singular fibres are either line pairs or double
lines.

 \begin{lem} In the above notation, $E$ has only singularities on the
singular fibres of the map $E\to C$. When the singular fibre is a line pair,
we have an $\rA{n}$ singularity at the point where the two components
meet (we include here the possibility $n=0$ when the point is a smooth point
of $E$). When the singular fibre is a double line, we have a
$\rD{n}$ singularity on the fibre (here we need to include the case $n=2$,
where we in fact have two $\rA1$ singularities, and $n=3$, where we
have an $\rA3$ singularity).
 \end{lem}

 \begin{pf} The proof is obvious, once the correct statement has been found.
The statement of this result in \cite{17} omits (for fibre a double line)
the cases $\rD{n}$ for $n>2$.
 \end{pf}

 \begin{lem}\label{lem_2.2} Suppose that $E\to C$ as above has $a_r$ fibres
which are line pairs with an $\rA{r}$ singularity and $b_s$ fibres which are
double lines with a $\rD{s}$ singularity (for $r\ge0$ and $s\ge2$), then 
 \[
 K_E^2=8(1-g)-\sum_{r\ge 0} a_r (r+1)-\sum_{s\ge 2} b_s s,
 \]
where $g$ denotes the genus of $C$.
 \end{lem}

This enables us to give a slick calculation of the Gromov--Witten invariants
when the base curve has genus $g>0$. In this case, it was shown in
\cite{17} that for a generic deformation of $X$, only finitely many fibres
from $E$ deform, and hence the Type~III contraction deforms to a Type~I
contraction. Thus Gromov--Witten numbers $n_1$ and $n_2$ may be defined as
in Section 1, and are nonnegative integers.

 \begin{prop}\label{prop_2.3} When $g>0$, we have 
 \[
 n_1=2\sum_{r\ge 0} a_r (r+1)+2\sum_{s\ge 2} b_s s \quad 
 \text{and} \quad n_2=2g-2.
 \]
 \end{prop}

 \begin{pf} We take a generic 1-parameter 
deformation of $X$, for which the
Type~III contraction deforms to a Type~I contraction. We therefore have a
diagram
 \[
 \renewcommand{\arraystretch}{1.3}
 \begin{matrix}
 \cX &\longrightarrow &\overline\cX \\
 \downarrow && \downarrow \\
 \De &=& \De
 \end{matrix}
 \]
where $\De\subset\C$ denotes a small disc. Since the singular locus of 
$\overline\cX $ consists only of curves of cDV singularities, we may again
apply \cite{8}, Theorem~11.10 to deduce the existence of a (smooth)
flopped fourfold $\cX'\to \overline\cX$. The induced family $\cX'\to\De$ is
given generically by flopping the fibres, and at $t=0$ it is easily checked
that $X_0'\iso X_0$; this operation is often called an {\em elementary
transformation} on the family. Identifying the groups $H^2 (X_t,\Z)\iso
H^2 (X_t',\Z)$ as before, this has the effect (at $t=0$) of sending $E$ to
$-E$ (cf.\ the discussion in \cite{5}, \S3.3). So if $E'$ denotes the class
in $H^2(X_t',\Z)$ corresponding to the class $E$ in
$H^2(X_t,\Z)$, we have $(E')^3=-E^3$. For $t\ne 0$, we just have a flop, and
so $(E')^3$ can be calculated from equation (2.1.1), namely
$(E')^3=E^3+n_1+8n_2$. Therefore, using Lemma~\ref{lem_2.2} 
 \[
n_1+8n_2=-2E^3=16(g-1)+2\sum_{r\ge0}a_r(r+1)+2\sum_{s\ge2}b_ss.
 \]

Similarly, we have $c_2(X')\cdot E'=- c_2(X)\cdot E$, and so from
equation (2.1.2) it follows that $2n_1+4n_2=2 c_2\cdot E$. An
easy calculation of the right-hand side then provides the second equation 
 \[
 2n_1+4n_2=8(g-1)+4\sum_{r\ge0}a_r(r+1)+4\sum_{s\ge2}b_ss.
 \]
 Solving for $n_1$ and $n_2$ from these two equations gives the desired
result. \end{pf}

 \begin{rem}\label{rem_2.4} This result remains true even when
$g=0$, although the slick proof given above is no longer valid. The formula
for $n_1$ is checked in Theorem~\ref{thm_2.5} by local deformation arguments
(for which the genus $g$ is irrelevant), showing that the contribution to
$n_1$ from a line pair fibre with $\rA{r}$ singularity is $2(r+1)$, and
from a double line fibre with $\rD{s}$ singularity is $2s$. Let $A\in
H_2(X,\Z)$ denote the class of a fibre of $E\to C$. Observe that any
pseudo\-holomorphic curve representing the numerical class $\eta$ will be a
component of a singular fibre of $E\to C$. Moreover, the components $l$ of a
singular fibre represent the same class in $H_2(X,\Z)$, and so in
particular twice this class is $A$. Thus the Aspinwall--Morrison formula (as
proved in \cite{15}) yields the contribution to the Gromov--Witten
invariants $\Phi_A (D,D,D)$ from double covers, purely in terms of $n_1$ and
$D\cdot A$. The difference may be regarded as the contribution to $\Phi_A
(D,D,D)$ from simple maps, and taking this to be $n_2 (D\cdot A)^3$
determines the number $n_2$ (in \S\ref{sec3}, we shall see how $n_2$ may be
determined directly from the moduli space of simple stable holomorphic
maps). If $g>0$, the above argument shows that this is in agreement with our
previous definition, and yields moreover the equality $n_2=2g-2$. The fact
that $n_2=-2$ when $g=0$ requires a rather more subtle argument involving
technical machinery -- see Theorem~\ref{thm_3.1}. I remark that the value
$n_2=-2$ is needed in physics, and that there is also a physics argument
justifying it (see \cite{4}, \S5.2 and \cite{5}, \S3.3) -- essentially, it
comes down to a statement about the A-model 3-point correlation functions.
In \S\ref{sec3} below, we give a rigorous mathematical proof of the
assertion.
 \end{rem}

 \begin{thm}\label{thm_2.5} The formula for $n_1$ in
Proposition~\ref{prop_2.3} is valid irrespective of the value of the genus
$g=g(C)$.
 \end{thm}

 \begin{pf} By making a holomorphic deformation of the complex structure on
an open neighbourhood $U$ in $X$ of the singular fibre $Z$ of $E\to C$, we
may calculate the contribution to $n_1$ from that singular fibre -- see
\cite{18}, (4.1). The deformation of complex structure is obtained as in
\cite{18} by considering the one dimensional family of Du Val singularities
in $\Xbar$, and deforming this family locally in a suitable neighbourhood
$\Ubar$ of the dissident point. Our assumption is that the family
$\Ubar\to\De$ has just an $\rA1$ singularity on $\Ubar_t$ for $t\ne 0$,
and we may assume also that $\Ubar\to \De$ is a good representative (in the
sense explained in \cite{18}). The open neighbourhood $U$ is then the blowup
of $\Ubar$ in the smooth curve of Du Val singularities (\cite{18}, p.\ 569).
The contribution to $n_1$ may be calculated locally, and will not change
when we make small holomorphic deformations of the complex structure on $U$,
which in turn corresponds to making small deformations to the family
$\Ubar\to\De$.

First we consider the case where the singular fibre $Z$ is a line pair --
from this, it will follow that the dissident singularity on $\Ubar$ is a
$\rcA{n}$ singularity with $n>1$, and that $\Ubar$ has a local analytic
equation of the form
 \[
x^2+y^2+z^{n+1}+tg(x,y,z,t)=0
 \]
in $\C^3\times\De$ (here $t$ is a local coordinate on $\De$, and $x=y=z=0$
the curve $C$ of singularities). For $t\ne 0$, we have an $\rA1$
surface singularity, which implies that $g$ must contain a term of the form
$t^r z^2$ for some $r\ge 0$. By an appropriate analytic change of
coordinates, we may then assume that $\Ubar$ has a local analytic equation of
the form
 \[
 x^2+y^2+z^{n+1}+t^{r+1}z^2+t h(x,y,z,t)=0,
 \]
where $h$ consists of terms which are at least cubic in $x,y,z$. By making a
small deformation of the family $\Ubar\to \De$, we may reduce to the case
$n=2$, that is, $\Ubar$ having local equation $x^2+y^2+z^3+t^{r+1}z^2+th=0$.
At this stage, we could in fact also drop the term $th$ (an easy check using
the versal deformation family of an $\rA2$ singularity), but this will
not be needed.

We now make a further small deformation to get $\Ubar_{\ep} \subset 
\C^3 \times \De $ given by a polynomial 
 \[
x^2+y^2+z^3+t^{r+1} z^2+\ep z^2+th=x^2+y^2+z^2 (z+t^{r+1}+\ep)+th \ .
 \]
This then has $r+1$ values of $t$ for which the singularity is an $\rA2$
singularity -- for other values of $t$, it is an $\rA1$ singularity. If we
blow up the singular locus of $\Ubar_{\ep}$, we therefore obtain a smooth
exceptional divisor for which $r+1$ of the fibres over $\De$ are line pairs.
By the argument of \cite{18}, (4.1), this splitting of the singular fibre
into $r+1$ line pair singular fibres of the simplest type can be achieved by
a local holomorphic deformation on a suitable open neighbourhood of the
fibre in the original threefold $X$.

It is however clear that a line pair coming from a dissident $\rcA2$
singularity of the above type contributes precisely two to the
Gromov--Witten number $n_1$ -- one for each line in the fibre. In terms of
equations, we have a local equation for $\Xbar$ of the form $x^2+y^2+z^3+w
z^2=0$; deforming this to say $x^2+y^2+z^3+w z^2+\ep w=0$, we get two simple
nodes, and hence two disjoint $(-1,-1)$-curves on the resolution.

The argument of \cite{18}, (4.1) shows that the Gromov--Witten number $n_1$
may be calculated purely from these local contributions, and so the total
contribution to $n_1$ from the line pair singular fibre of $E$ with 
$\rA{r}$ singularity is indeed $2(r+1)$, as claimed.

For the case of the singular fibre $Z$ of $E$ being a double line, the
dissident singularity must be $\rcE6$, $\rcE7$, $\rcE8$, or $\rcD{n}$ for
$n\ge4$. Thus $\Ubar$ has a local analytic equation of the form
$f(x,y,z)+tg(x,y,z,t)$ in $\C^3\times\De$ for $f$ a polynomial of the
appropriate type ($t$ a local coordinate on $\De$, and $x=y=z=0$ the curve
of singularities). To simplify matters, we may deform $f$ to a polynomial
defining a $\rD4$ singularity, and hence make a small deformation of the
family to one in which the dissident singularity is of type $\rcD4$. We then
have a local analytic equation of the form
 \[
x^2+y^2 z+z^3+t g(x,y,z,t)=0.
 \]

For $t\ne 0$, we have an $\rA1$ singularity, and so the terms of $g$ must be
at least quadratic in $x,y,z$. Moreover, by changing the $x$-coordinate, we may
take the equation to be of the form 
 \[
x^2+y^2 z+z^3+t^a y^2+t^b yz+t^c z^2+t h(x,y,z,t)=0,
 \]
with $a,b,c$ positive, and where the terms of $h$ are at least cubic in $x,y,z$. 
The fact that the blowup $U$ of $\Ubar$ in $C$ is smooth is easily checked to
imply that $a=1$. Since
 \[
ty^2+2 t^b yz=t(y+t^{b-1}z)^2-t^{2b-1}z^2,
 \]
 we have an obvious change of $y$-coordinate which brings the
equation into the form
 \[
x^2+y^2 z+z^3 +t y^2+t^r z^2+t h_1
(x,y,z,t)=0,
 \]
 where $r=\min \{ c, 2b-1 \}$ and $h_1$ has the same
property as $h$.

When we blow up $\Ubar$ along the curve $x=y=z=0$, we obtain an exceptional locus
$E$ with a double fibre over $t=0$, on which we have a $\rD{r+1}$ singularity (including the case $r=1$ of two $\rA1$ singularities,
and $r=2$ of an $\rA3$ singularity). Moreover, this was also true of our
original family, since the small deformation of $f$ we made did not affect the
local equation of the exceptional locus.

Moreover, by adding a term $\ep_1 y^2+\ep_2 z^2$, we may
deform our previous equation to one of the form 
 \[
 x^2+y^2(z+t+\ep_1)+z^2(z+t^r+\ep_2)+th_1(x,y,z,t)=0.
 \]
 When $ t+\ep_1=0$, we have an $\rA3$ singularity, and when $t^r+\ep_2=0$,
an $\rA2$ singularity. Moreover, when we blow up the singular locus of this
deformed family, the resulting exceptional divisor is smooth and has line
pair fibres for these $r+1$ values of $t$. Thus, as seen above, the
contribution to $n_1$ from the original singular fibre (a double line with a
$\rD{r+1}$ singularity) is $2(r+1)$ as claimed.
 \end{pf}

 \section{Calculation of $n_2$ for Type~III contractions}\label{sec3}
 Let $\fiemapX $ be a Type~III contraction on a Calabi--Yau threefold $X$,
which contracts a divisor $E$ to a (smooth) curve $C$ of genus $g$. When
$g>0$, it was proved in Proposition~\ref{prop_2.3} that the Gromov--Witten
number $n_2$ (defined for arbitrary genus via Remark~\ref{rem_2.4}) is
$2g-2$. The purpose of this Section is to extend this result to include the
case $g=0$ ($C$ is isomorphic to $\PP^1$), and to prove $n_2=2g-2$ in
general. 

Arguing as in \cite{18}, it is clear that the desired result is a local one,
depending only on a neighbourhood of the exceptional divisor $E$. As remarked
in \S\ref{sec2}, we may then always reduce down to the case that the generic
fibre of $E\to C$ is irreducible. If all the fibres of $E\to C$ are smooth (so
$E$ is a $\PP^1$-bundle over $C$), the fact that $n_2=2g-2$ was proved in
Proposition 5.7 of \cite{11}, using a cobordism argument. This latter result
was extended by Ruan in \cite{13}, Proposition~2.10, using the theory of
moduli spaces of stable maps and the virtual neighbourhood technique
(cf.~\cite{2,9}). If the singular fibres of $E\to C$ are line pairs, Ruan's
result applies directly. We prove below that the linearized Cauchy--Riemann
operator has constant corank for the stable (unmarked) rational curves given
by the fibres of $E$ over $C$, and hence by Ruan's result that there is an
obstruction bundle $\cH$ on $C$, with $n_2$ determined by the Euler class of
$\cH$. By Dolbeault cohomology, there is a natural identification of
$\cH$ with the cotangent bundle $T_C^*$ on $C$, and hence the formula for $n_2$
follows. We note however that for Ruan's result to hold, we do not need an 
integrable almost complex structure on $X$. Provided we have a natural
identification between the cokernel of the linearized Cauchy--Riemann operator
and the cotangent space at the corresponding point of $C$, we can still deduce
that $n_2=2g-2$. In the general case of a Type~III contraction which has
double fibres, we show below that we can make a small local deformation of the
almost complex structure on $X$ so that $E$ deforms to a family of
pseudo\-holomorphic rational curves over $C$ with at worst line pair singular
fibres, and for which Ruan's method applies.

 \begin{thm}\label{thm_3.1} For any Type~III contraction $\fiemapX $, 
the Gromov--Witten number $n_2=2g-2$.
 \end{thm}

 \begin{pf} We saw above that we may assume that the generic fibre of $E\to
C$ is irreducible. Furthermore, we initially assume also that the singular
fibres are all line pairs, and later reduce the general case to this one.

We let $J$ denote the almost complex structure on $X$, which we know is
integrable (at least in a neighbourhood of $E$), and tamed by a symplectic
form $\omega$. Let $A\in H_2(X,\Z)$ be the class of a fibre of $E\to C$.
Adopting the notation from \cite{13}, we consider the moduli space 
$\cMbar_A(X,J)=\cMbar_A(X,0,0,J)$ of stable unmarked rational
holomorphic maps, a compactification of the space of (rigidified)
pseudo\-holomorphic maps $\C\PP^1\to X$, representing the class
$A$. The theory of stable maps, as explained in Section 3 of \cite{13}, goes
through for unmarked stable maps, by taking each component of the domain as
a bubble component, and adding marked points (in addition to the double
points) as in \cite{13} in order to stabilize the components (thus taking a
local slice of the automorphism group).

In the case that all the singular fibres of $E\to C$ are line pairs,
$\cMbar_A (X,J)$ has two components, one corresponding to simple maps and
the other to double covers. It is now a simple application of Gromov
compactness to see that these two components are disjoint, since a sequence
of double cover maps cannot converge to a simple map. A similar argument
will show that for all almost complex structures $J_t$ in some neighbourhood
of $J=J_0$, the moduli space $\cMbar_A (X,J_t)$ will consist of two
disjoint components, one corresponding to the simple maps and the other to
the double covers.

Since any stable unmarked rational holomorphic map must be an embedding, it
is clear that the component $\cMbar'_A (X,J)$ corresponding to the simple
maps can be identified naturally with the smooth base curve $C$, and that
for all almost complex structures in some neighbourhood of $J=J_0$, the
moduli space $\cMbar'_A (X,J_t)$ of simple unmarked stable holomorphic maps
is compact. The Gromov--Witten invariant $n_2$ that we seek can then be
defined via Ruan's virtual neighbourhood invariant $\mu_{\cS}$, and may be
evaluated on $(X,J)$ by using \cite{13}, Proposition~2.10.

Let us now go into more details of this. We consider $C^{\infty}$ stable
maps $f\in\Bbar_A (X)=\Bbar_A (X,0,0)$ in the sense of \cite{13}, Definition
3.1, where Ruan shows later in the same Section that the naturally
stratified space $\Bbar_A (X)$ satisfies a property which he calls {\em
virtual neighbourhood technique admissable} or {\em VNA}, and as he says, for
the purposes of the virtual neighbourhood construction, behaves as if it were a
Banach $V$-manifold. Since any simple marked holomorphic stable map $f$ in
$\cMbar'_A (X,J)$ is forced to be an embedding, we may restrict our
attention to $C^{\infty}$ stable maps whose domain $\Si$ comprises at most
two $\PP^1$s. We stratify $\Bbar_A (X)$ according to the combinatorial type
$D$ of the domain $\Si$. Thus any $f\in \cMbar'_A (X,J)$ belongs to one of
two strata of $\Bbar_A (X)$.

In general, for $k$-pointed $C^{\infty}$ stable maps of genus $g$, Ruan shows
that for any given combinatorial type $D$, the substratum $ B_D (X,g,k)$ is
a Hausdorff Frechet V-manifold (\cite{13}, Proposition 3.6). As mentioned
above, he needs to add extra marked points in order to stabilize the
nonstable components of the domain $\Si$, thus taking a local slice of the
action of the automorphism group on the unstable marked components of $\Si$.
Moreover, the tangent space $T_f B_D (X,g,k)$ is identified with $\Om^0 (f^*
T_X)$, as defined in his equation \cite{13}, (3.29). The tangent space $T_f
\Bbar_A (X,g,k)$ can then be defined as $T_f B_D (X,g,k) \times \C_f $,
where $\C_f$ is the space of gluing parameters (see \cite{13}, equation
before (3.67)).

In our case, however, things are a bit simpler. Given $f\in \cMbar'_A
(X,J)$ with domain $\Si$ consisting of two $\PP^1$s, the tangent space $T_f
\Bbar_A (X)$ is of the form $\Om^0 (f^* T_X) \times \C$, and we have a
neighbourhood $\wave U_f$ of $f$ in $\Bbar_A (X)$ defined by
\cite{13}, (3.43), consisting of stable maps $\overline f^{v,w}$ parametrized
locally by
 \[
 \bigl\{w\in\Om^0(f^*T_X)\ ;\ \|w\|_{C^1}<\ep'\bigr\}
 \]
(corresponding to deformations within the stratum $B_D(X)$), and by
$v\in\C_f^\ep $ (an $\ep$-ball in $\C_f=\C$ giving the gluing parameter at
the double point). This then corresponds to the above decomposition of
$T_f\Bbar_A (X)$ into two factors. On the first factor, the linearization
$D_f\delbar_J $ of the Cauchy--Riemann operator restricts to
 \[
\delbar_{J,f}\colon\Om^0(f^*T_X)\to\Om^{0,1}(f^*T_X)
 \]
in the notation of \cite{13}. The index of this operator may be calculated
using Riemann--Roch on each component of $\Si$ (cf.~the proof of Lemma
3.16 in \cite{13}, suitably modified to take account of the extra marked
points), and is seen to be $-2$. 

Let us now consider the stable maps $f^v=f^{v,0}$ for
$v\in\C_f^\ep\setminus\{0\}$. These are stable maps $\C\PP^1\to X$ which
differ from $f$ only in small discs around the double point, and in this
sense are approximately holomorphic. Set $v=r e^{i\theta}$; then the
gluing to get $f^v\colon \Si^v\to X$ is only performed in discs around the
double point of radius $2r^2 /\rho$ in the two components ($\rho$ a suitable
constant). It can then be checked for any $2<p<4$ that $\| \delbar_J (f^v)\|_{L^p_1} \le Cr^{4/p}$ (see \cite{13} Lemma 3.23, and
\cite{10} Lemma A.4.3), from which it follows that the linearization 
 \[
L_A=D_f \delbar_J 
 \]
of the Cauchy--Riemann operator should be taken as zero on the factor $\C_f$
in $T_f \Bbar_A (X)$. Thus we deduce that the index of $L_A$ is zero, and
that $\coker L_A$ is same as the cokernel of $\delbar_{J,f}\colon \Om^0 (f^*
T_X)\to \Om^{0,1} (f^* T_X)$, which by Dolbeault cohomology may be
identified as
 \[
 H^1(f^*T_X)=H^1(Z,T_X{}\rest Z),
 \]
where $Z$ is the fibre of $E\to C$ (over a point $x\in C$) corresponding to
the image of $f$.

We note that these are exactly the same results as are obtained in the smooth
case, when $\Si$ consists of a single $\PP^1$. Here, we need to add three
marked points to stabilize $\Si$, and Riemann--Roch then gives immediately
that the index of $L_A$ is zero.

Observe that $Z$ is a complete intersection in $X$, and so for our purposes
is as good as a smooth curve. Via the obvious exact sequence,
$H^1(T_X{}\rest Z)$ may be naturally identified with $H^1(N_{Z/X})$, which
in turn may be naturally identified with $H^0 (N_{Z/X})^*$ (since
$K_Z=\bigwedge^2 N_{Z/X}$, we have a perfect pairing $H^0 (N_{Z/X}) \times
H^1(N_{Z/X})\to H^1(K_Z)\iso \C$). Observing that $H^0 (N_{Z/X})=H^0
(\Oh_Z\oplus\Oh_Z (E))\iso \C$, we know that $\coker L_A $ has complex
dimension one and is naturally identified with $T^*_{C,x}$, the dual of the
tangent space at $x$ to the Hilbert scheme component $C$. This we have seen
is true for all $f\in\cMbar'_A(X,J)$.

We now apply \cite{13}, Proposition~2.10, (2) to our set-up, where
$C=\cMbar'_A (X,J)=\cM_{\cS}=\cS^{-1}(0)$ for $\cS$ the Cauchy--Riemann
section of $\overline \cF_A(X)$ (as constructed in \cite{13}, \S3) over a
suitable neighbourhood of $\cM_{\cS}$ in $\Bbar_A(X)$. The above calculations
verify that the conditions of Proposition~2.10, (2) are satisfied, with
$\ind(L_A)=0$, $\dim(\coker L_A)=2$ and $\dim (\cM_{\cS})=2$. Moreover, we
deduce that the obstruction bundle $\cH$ on $\cM_{\cS}$ is just the cotangent
bundle $T_C^*$ on $C$.

The Gromov--Witten number $n_2$ may then be defined to be $\mu_{\cS}(1)$. 
It follows from the basic Theorem~4.2 from \cite{13} that this is
independent of any choice of tamed almost complex structure and is a
symplectic deformation invariant. Thus by considering a small deformation of
the almost complex structure and using \cite{13}, Proposition~2.10, (1), it
is the invariant $n_2$ that we seek. Applying Ruan's crucial
Proposition~2.10, (2), the invariant can be expressed as
 \[
 \mu_{\cS}(1)=\int_{\cM'_A (X,J)} e(T_C^*),
 \]
from which it follows that $n_2=2g-2$ as claimed.

The general case (where $E\to C$ also has double fibres) can now be reduced
to the case considered above. Suppose we have a point $Q\in C$ for which the
corresponding fibre is a double line. We choose an open disc $\De \subset C$
with centre $Q$, and a neighbourhood $U$ of $Z$ in $X$, with $U$ fibred over
$\De$, the fibre $U_0$ over $Q$ containing the fibre $Z$. Letting $\Ubar\to
\De$ denote the image of $U$ under $\fie$, a family of surface Du Val
singularities, we make a small deformation $\overline\cU\to \De'$ of $\Ubar$, 
as in the proof of Theorem~\ref{thm_2.5} of this paper, and in this way
obtain a holomorphic deformation $\cU\to \De'$ of $U$ under which
$E_0=E\rest{\De}$ deforms to a family of surfaces $E_t$ ($t\in \De'$), all
fibred over $\De$, and with at worst line pair singular fibres for $t\ne 0$.
Considering $\overline\cU\to \De \times \De'$ as a two parameter deformation
of the surface singularity $ \Ubar_0$, we may take a good representative and
apply Ehresmann's fibration theorem (with boundary) to the corresponding
resolution $\cU\to\De\times\De'$ (cf.\ \cite{18}, proof of Lemma 4.1).
In this way, we may assume that $\cU\to \De \times \De'$ is differentiably
trivial over the base. In particular, the family $\cU\to\De'$ is also
differentiably trivial, and hence determines a holomorphic deformation of
the complex structure on a fixed neighbourhood $U$ of $Z$, where $U\to \De$
is also differentiably trivial.

We perform this procedure for each singular fibre $Z_1,\dots,Z_N $ of
$E\to C$, obtaining, for each $i$, an open neighbourhood $U_i$ of $Z_i$
fibred over $\De_i \subset C$, and a holomorphic complex structure $J_i$ on
$U_i$ with the properties explained above (of course, if $Z_i$ is a line
pair, we may take $J_i$ to be the original complex structure $J$). Let
$\frac{1}{2}\overline\De_i$ denote the closed subdisc of $\De_i$ with half the
radius, $C^*=C\setminus\bigcup_{i=1}^N\frac{1}{2}\overline\De_i$, and
$E^*=E\rest{C^*}\to C^*$ the corresponding open subset of $E$. We then take
a tubular neighbourhood $U^*\to C^*$ of $E^*\to C^*$, equipped with the
original complex structure $J$. By taking deformations to be sufficiently
small and shrinking radii of tubular neighbourhoods if necessary, all these
different complex structures may be patched together in a $C^{\infty}$ way
(tamed by the symplectic form) over the overlaps in $C$. In this manner, we
obtain an open neighbourhood $W$ of $E$ in $X$, and a tamed almost complex
structure $J'$ on $W$, which is a small deformation of the original complex
structure $J$ and which satisfies the following properties:

 \begin{enumerate}
 \renewcommand{\labelenumi}{(\alph{enumi})}
 \item Each singular fibre $Z_i$ of $E\to C$ has an open neighbourhood $U_i
\subset W$ fibred over $\De_i \subset C$ with $J'$ inducing an integrable
complex structure on each fibre (thus $U_i\to \De_i$ is a $C^{\infty}$
family of holomorphic surface neighbourhoods).

 \item The almost complex structure $J'$ is integrable in a smaller 
neighbourhood $U'_i \subset U_i$ of each singular fibre, with the
corresponding family $U'_i\to \De'_i$ being holomorphic.

 \item On the complement of $\bigcup U_i$ in $W$, the almost complex
structure $J'$ coincides with the original complex structure $J$.

 \item $E$ deforms to a $C^{\infty}$ family of pseudo\-holomorphic rational
curves $E'\to C$ in $(W,J')$, with generic fibre $\C\PP^1$ and the only
singular fibres being line pairs. Moreover, we may assume that any such
singular fibre is contained in one of the above open sets $U'_i$.
 \end{enumerate}

Of course, we may now patch $J'$ on $W$ with the original complex structure $J$
on $X$ to get a global tamed almost complex structure on $X$, which we shall also
denote by $J'$. Provided we have taken our deformations sufficiently small, the
standard argument via Gromov compactness ensures that any pseudo\-holomorphic
stable map representing the class $A$ has image contained in a fibre of $E'\to
C$.

The theory of \cite{13} applies equally well to almost complex structures, and hence
to our almost complex manifold $X'$ with complex structure $J'$. Clearly, all
the calculations remain unchanged for stable maps whose image (a fibre of 
$E'\to C$) has a neighbourhood on which $J'$ is integrable, and in particular
this includes all the singular fibres. Suppose therefore that $f\colon \C\PP^1\to 
X'$ is a pseudo\-holomorphic rational curve whose image $Z$ is contained in an
overlap $U_i \setminus U'_i$ (where $J'$ may be nonintegrable). The linearized 
Cauchy--Riemann operator $L_A$ still has index zero, since by the argument
of \cite{10}, p.~24, the calculation via Riemann--Roch continues to give the
correct value. We therefore need to show that $\coker L_A$ is still
identified naturally as $T^*_{C,x}$, and hence that the obstruction bundle 
is $\cH=T^*_C$ as before.

Setting $U=U_i$ and $\De=\De_i$, we know that $U\to \De$ is locally
(around the image $Z$ of $f$) a $C^{\infty}$ family of holomorphic surface
neighbourhoods. Moreover, the linearized 
Cauchy--Riemann operator $L_A = D_f\colon C^{\infty} (f^*T_U)\to \Om^{0,1} 
(f^*T_U)$ fits into the
following commutative diagram (with exact rows)
 \[
 \renewcommand{\arraystretch}{1.3}
 \begin{matrix}
 0&\to&C^{\infty}(f^*T_{U/\De})&\to
 &C^{\infty}(f^*T_U)&\to&C^{\infty}(g^*T_{\De})&\to&0 \\
 &&\Bigm\downarrow\delbar_{f}&&\Bigm\downarrow D_f&&\Bigm\downarrow&& \\
 0&\to&\Om^{0,1}(f^*T_{U/\De})&\to
 &\Om^{0,1}(f^*T_U)&\to&\Om^{0,1}(g^*T_{\De})&\to&0
 \end{matrix}
 \]
where $g$ is the constant map on $\C\PP^1$ with image the point $x\in \De$,
and where the fibre of $E'$ over $x$ is $Z$. Let us denote by $U_x$ the
corresponding holomorphic surface neighbourhood, the fibre of $U$ over $x$. The
cokernel of 
 \[
 \delbar_{f}\colon C^{\infty}(f^*T_{U/\De})\to\Om^{0,1}(f^*T_{U/\De})
 \]
is then naturally identified via Dolbeault cohomology with
$H^1(T_{U_x}{}\rest Z)\iso H^1(N_{Z/U_x})$. This latter space is in turn
naturally identified with $H^1(N_f)\iso H^0(N_f)^*\iso T^*_{C,x}$.

I claim now that $J'$ may be found as above for which $\coker L_A$ has the
correct dimension (namely real dimension two) for all fibres of $E'\to C$.
Since $L_A$ has index zero and $\ker L_A$ has dimension at least two, we 
need to show that that the dimension of $\coker L_A$ is not more than two.  
This follows by a Gromov compactness argument. Suppose that the dimension is 
too big for some fibre of $E'\to C$, however close we take $J'$ to $J$. We
can then find sequences of almost complex structures $J'_\nu$ (with the
properties (a)--(d) described above) converging to $J=J_0$, and
pseudo\-holomorphic rational curves $f_\nu\colon\C\PP^1\to (X, J'_\nu)$ at
which $\coker L_A$ has real dimension $>2$. By construction, the image of
such a map is not contained in any $U'_i$ (since $J'_\nu$ would then be
integrable on some neighbourhood of the image, and then we know that $\coker
L_A$ has the correct dimension). Thus the image of $f_\nu$ has nontrivial
intersection with the compact set $X \setminus \bigcup U'_i$. By Gromov
compactness, the $f_\nu$ may be assumed to converge to a pseudo\-holomorphic
rational curve on $(X,J)$ whose image is not contained in any $U'_i$. This is
therefore just an embedding $f\colon\C\PP^1\to (X, J)$ of some smooth fibre
of $E\to C$, at which we know that $\coker L_A$ has real dimension
precisely two; this then is a contradiction. A similar argument, via Gromov
compactness, then yields the fact that $J'$ may be found as above such that
the linear map $\coker (\delbar_{f})\to \coker (D_f)$ is an isomorphism 
for all smooth fibres of $E'\to C$, since this is
true for all the smooth fibres of $E\to C$ on $(X,J)$.

For such a $J'$, we deduce that $\coker L_A$ is naturally identified with 
$T^*_{C,x}$ for all fibres, and hence the obstruction bundle identified as 
$T^*_C$. The previous argument may then be applied directly with the almost
complex structure $J'$, showing that the symplectic invariant $n_2$ is
$2g-2$ in general. The proof of Theorem~\ref{thm_3.1} is now complete.
\end{pf}

 \section{An application to symplectic deformations of Calabi--Yaus}
\label{sec4}
 If $X$ is a Calabi--Yau threefold which is general in moduli, we know that
any codimension one face of its nef cone $\Kbar (X)$ (not contained in the
cubic cone $W^*$) corresponds to a primitive birational contraction
$\fiemapX $ of Type~I, II or $\text{III}_0$, where Type~$\text{III}_0$
denotes a Type~III contraction for which the genus of the curve $C$ of
singularities on $\Xbar$ is zero.

In \cite{18}, we studied Calabi--Yau threefolds which are symplectic
deformations of each other. One of the results proved there (Theorem~2) said
that if $X_1$ and $X_2$ are Calabi--Yau threefolds, general in their complex
moduli, which are symplectic deformations of each other, then their K\"ahler
cones are the same. The proof of this essentially came down to showing that
certain Gromov--Witten invariants associated to exceptional classes were
nonzero. Using the much more precise information obtained in this paper, we
are able to make a stronger statement.

 \begin{cor}\label{cor4.1} With the notation as above, any codimension one
face (not contained in $W^*$) of $\Kbar (X_1)=\Kbar (X_2)$ has the same
contraction type (Type~I, II or $\text{\sl III}_0$) on $X_1$ as on
$X_2$.
 \end{cor}

 \begin{pf} The fact that Type~II faces correspond is easy, since for $D$ in
the interior of such a face, the quadratic form $q(L)=D\cdot L^2$ is
degenerate, which is not the case for $D$ in the interior of a Type~I or
Type~$\text{III}_0$ face. Stating it another way, if we consider the Hessian
form associated to the topological cubic form $\mu$, then $h$ is a form of
degree $\rho=b_2$ which has a linear factor corresponding to each Type~II
face. Thus the condition that a face is of Type~II is topologically
determined.

The result will therefore follow if we can show that a face of the nef cone
which is Type~I for one of the Calabi--Yau threefolds is not of
Type~$\text{III}_0$ for the other. However, for a Type~I face, we saw in
\S\ref{sec1} that $n_d$ is always nonnegative; for a Type~$\text{III}_0$
face, we proved in Theorem~\ref{thm_3.1} that $n_2=-2$. Since Gromov--Witten
invariants are invariant under symplectic deformations, the result is proved.
 \end{pf}

 \begin{rem} It is still an open question whether there exist examples of
Calabi--Yau threefolds $X_1$ and $X_2$ which are symplectic deformations of
each other but not in the same algebraic family.
 \end{rem}

\noindent Mathematics Subject Classification (1991):

14J10, 14J15, 14J30, 32J17, 32J27, 53C15, 53C23, 57R15, 58F05 

\bigskip
\noindent
P.M.H. Wilson \\
Department of Pure Mathematics, University of Cambridge, \\
16 Mill Lane, Cambridge CB2 1SB, UK \\
e-mail: pmhw@dpmms.cam.ac.uk

 \end{document}